\documentclass{article}
\usepackage{amsmath,amsthm,amssymb,amsfonts,amscd,pstricks}
\begin{document}

{\bf \Large Redundancies : An Omission in Probability \\ \\ Theory ?} \\ \\

{\bf Elem\'{e}r E Rosinger} \\
Department of Mathematics \\
and Applied Mathematics \\
University of Pretoria \\
Pretoria \\
0002 South Africa \\
eerosinger@hotmail.com \\ \\

{\bf Abstract} \\

It is shown that the standard Kolmogorov model for probability spaces cannot in general allow
the elimination but of only a small amount of probabilistic redundancy. This issue, a purely
theoretical weakness, not necessarily related to empirical reality, appears not to have
received enough attention in foundational studies of Probability Theory. \\ \\

{\bf 1. The Need for a Match} \\

Several recent works have brought to attention the events during the first decades of the 20th
century which culminated in 1933 with Kolmogorov's foundation of modern Probability Theory in
his {\it Grundbegriffe}, see Shafer \& Vovk [1], Vovk \& Shafer, von Plato, or Bilova et.al. \\

Ever since Probability Theory first started to emerge in the work of Jacob Bernoulli early in
the 18th century, there has been an interest in the extent to which such a theory may in fact
match with the empirical reality it is supposed to model. \\

The axiom regarding "total probability", namely \\

$~~~~~ P ( A \cup B ) ~=~ P ( A ) + P ( B ) $ \\

for two events $A$ and $B$ which cannot happen simultaneously, or the axiom for "compound
probability", namely \\

$~~~~~~ P ( A \cap B ) ~=~ P ( A )\, P ( B \,|\, A ) $ \\

for two arbitrary events did not lead to significant problems when related to empirical
reality. \\

What on the other hand led to such problems was the so called {\it Cournot Principle} which in
one of its formulations, Cournot, stated that :

\begin{quote}

{\it The physically impossible event is therefore the one that has infinitely small
probability, and only this remark gives substance - objective and phenomenal value - to the
theory of mathematical probability.}

\end{quote}

This idea, however, did not originate with Cournot. Jacob Bernoulli in his 1713 {\it Ars
Conjectandi} proved what would later be called the Law of Large Numbers which in simple terms
says that in a large enough sequence of independent trials of the same event there is a high
probability that the {\it frequency} of the event will be close to its {\it probability}. And
Bernoulli commented that we can consider that high probability as a certainty, and thus use
frequency as an estimate of probability. \\ \\

{\bf 2. Falsifiability} \\

Only recently in the history of science, and owing to Karl Popper, did the idea clearly arise
and got expressed of a proper test of the possible match, or otherwise, between a given
scientific theory and the empirical reality it is supposed to model. And the test of such a
match which Popper suggested, and called "falsifiability", consists of a set of fundamental
experiments clearly specified by the respective theory, experiments which can effectively be
performed within the empirical reality. And in case any of them may happen to fail, then the
theory as a whole is considered inadequate. \\

Obviously, two conditions should be met by any such set of fundamental experiments :

\begin{itemize}

\item the set of required experiments must be finite, and preferably small,

\item each such experiment should only require a finite time to perform, and preferably, a
short enough time at that.

\end{itemize}

When for instance, Einstein came up in 1915 with his theory of General Relativity, he himself
suggest several such fundamental experiments. One of them was the celebrated bending of light
rays near to a massive mass. And early in 1919, during a Sun eclipse, Arthur Eddington managed
to show that General Relativity passed that test, since the amount of bending predicted by
that theory was in good correspondence with the empirical bending measured. \\

Needless to say, the various theories of Physics, including General Relativity, have so far
been required to satisfy experiments which happen to conform to the above two
conditions. \\

In this regard, one of the major deficiencies of String Theory is precisely the fact that, so
far, it cannot offer any experiment which could effectively be performed and respect the above
two conditions. \\ \\

{\bf 3. The Special Situation of Probability Theory} \\

As far as the mentioned Cournot Principle is concerned, it is quite obvious that any empirical
experiments which would be able to test it in the sense of Popper's falsifiability are duty
bound to fail the second above condition. Indeed, there is simply no a priori way of knowing
how large a finite sequence of independent experiments of the same event should be, in order
to give the possibility for obtaining a sufficiently reliable acceptance or rejection of that
principle. \\

In this way, Probability Theory is in the special position to be unable to either confirm or
falsify itself on effective empirical grounds.
And this comes from the simple fact that aiming to establish a close enough connection between
{\it frequency} and {\it probability} cannot lead to empirical tests which satisfy the above
second condition, since any such close enough connection can only stand or fail at the {\it
infinite limit}. \\

Needless to say, even if not clearly and explicitly enough, there has been an awareness of
this problem, and consequently, a variety of solutions has been proposed, as seen for instance,
in Shafer \& Vovk [1]. \\

Among the more recent major ideas in this regard, one can mention the following two :

\begin{itemize}

\item Relating probability to {\it complexity} in the sense of Kolmogorov.

\item Relating probability to {\it games}, Shafer \& Vovk [2].

\end{itemize}

The first idea, already about half a century old, has proved to be particularly seminal, and
not only in Probability Theory, but also in several major branches of Mathematics, as seen for
instance in Chaitin. \\

The second idea, a more recent one, may at first appear as risking to lead well beyond what
would conventionally be understood by randomness or probability. Indeed, Kolmogorov complexity
is defined by the performance of the shortest possible program on a universal Turing machine.
On the other hand, a game, as suggested for instance in Shafer \& Vovk [2], Vovk \& Shafer,
comprises more than one {\it conscious} player. Thus one may become concerned that the
processes which such a game can generate may in fact go beyond what one usually conceives of
as being random or probabilistic. \\ \\

\newpage

{\bf 4. Redundancies : An Omission ?} \\

It appears that, historically, the concerns related to Probability Theory have been focused
upon its possible match with empirical reality, and dealt mainly, if not in fact, mostly with
the issue of the relation between frequency and probability, a relation which as mentioned
above, proves not to be within the realms of customary empirical tests. \\

And during this long ongoing process, it appears that an important {\it internal purely
theoretical weakness} of Probability Theory - a weakness which is clearly there, even if no
relation to empirical reality is involved - was not given attention. Namely, let us consider
an arbitrary probability space \\

(1) $~~~~~~ ( \Omega, {\cal F}, P ) $ \\

where ${\cal F}$ is a $\sigma$-algebra of subsets of $\Omega$, while$P : {\cal F}
\longrightarrow [ 0, 1 ]$ is the respective probability. Further, let us assume, as customary,
that \\

(1a) $~~~~~~ \{~ \omega ~\} \in {\cal F},~~~ \mbox{for}~~ \omega \in \Omega $ \\

There are now two cases :

\begin{itemize}

\item the set $\Omega$ is finite or countable,

\item the set $\Omega$ is uncountable.

\end{itemize}

In the first case, let us take \\

(2) $~~~~~ \Omega_1 ~=~ \{~ \omega \in \Omega ~~|~~ P ( \{~ \omega ~\} ) > 0 ~\} $ \\

in other words, we discard from $\Omega$ all the points $\omega$ which have probability zero,
thus are {\it redundant probabilistically}. In this way $\Omega_1$ represents the set of
probabilistically non-redundant points. Then clearly \\

$~~~~~~ \Omega_1 \in {\cal F} $ \\

since $\Omega_1 = \Omega \setminus \{~ \omega \in \Omega ~~|~~ P ( \{~ \omega ~\} ) = 0 ~\}$,
and the second set in the right hand term is at most countable, thus in view of (1a), it
belongs to ${\cal F}$. \\

Further, corresponding to $\Omega_1$, let us take \\

(3) $~~~~~~ {\cal F}_1 ~=~ \{~ A \cap \Omega_1 ~~|~~ A \in {\cal F} ~\} $ \\

thus ${\cal F}_1 \subseteq {\cal F}$, and $P ( B ) > 0$, for $B \in {\cal F}_1,~ B \neq
\phi$. \\

Finally, let us take \\

(4) $~~~~~~ P_1 : {\cal F}_1 ~\longrightarrow~ [ 0, 1 ]~~~
                                \mbox{defined by}~~ P_1 = P\,|\,_{{\cal F}_1} $ \\

In this way we obtain the {\it probabilistic redundancy free} probability space \\

(5) $~~~~~~ ( \Omega_1, {\cal F}_1, P_1 ) $ \\

And now one can note that much of the probabilistic phenomena of interest are {\it the same}
on the probability spaces $( \Omega, {\cal F}, P )$ and $( \Omega_1, {\cal F}_1, P_1 )$. And
needless to say, it is more convenient to work on the latter, precisely owing to the lack of
redundancies. \\

On the other hand, in the second case above, namely, when $\Omega$ is uncountable, one {\it
cannot} in general apply the above procedure in order to get rid of {\it probabilistic
redundancies}. For instance, let us take the simplest usual case when \\

(6) $~~~~~~ \Omega ~=~ [ 0, 1 ] ~\subset~ \mathbb{R} $ \\

while ${\cal F}$ and $P$ are, respectively, the set of Lebesgue measurable subsets of $[ 0,
1]$ and the Lebesgue measure on $[ 0, 1 ]$. Indeed, in this case we obtain \\

(7) $~~~~~~ \Omega_1 ~=~ \phi $ \\

since every point $x \in [ 0, 1 ]$ has Lebesgue measure zero. Thus we simply cannot eliminate
{\it all probabilistically redundant} points, since we do no longer remain with a nontrivial
probability space. \\

This failure related to the elimination of probabilistic redundancy in the Kolmogorov model in
one of the main reasons in the technical difficulties encountered in the study of stochastic
processes with continuous time. \\ \\

\end{document}